*Chapter*

# THE LINEARIZATION METHODS AS A BASIS TO DERIVE THE RELAXATION AND THE SHOOTING METHODS

### István Faragó[*], DSc

Department of Differential Equations, Institute of Mathematics, Budapest University of Technology and Economics, & MTA-ELTE Research Group, Budapest, Hungary

### Stefan M. Filipov, PhD

Department of Computer Science, University of Chemical Technology and Metallurgy, Sofia, Bulgaria

## ABSTRACT

This chapter investigates numerical solution of nonlinear two-point boundary value problems. It establishes a connection between three important, seemingly unrelated, classes of iterative methods, namely: the linearization methods, the relaxation methods (finite difference methods), and the shooting methods. It has recently been demonstrated that using finite differences to discretize the sequence of linear problems obtained by quasi-linearization, Picard linearization, or constant-slope linearization, leads to the usual iteration formula of the respective relaxation method. Thus, the linearization methods can be used as a basis to derive the relaxation methods. In this work we demonstrate that the shooting methods can be derived from the linearization methods, too. We show that relaxing a shooting trajectory, i.e. an initial value problem solution, is in fact a projection transformation. The obtained function, called projection trajectory, can be used to correct the initial condition. Using the new initial condition, we can find a new shooting trajectory, and so on. The described procedure is called shooting-projection iteration (SPI). We show that using the quasi-linearization equation to relax (project) the shooting trajectory leads to the usual shooting by Newton method, the constant-slope linearization leads to the usual shooting by constant-slope method, while the Picard linearization leads to the recently proposed shooting-projection method. Therefore, the latter method can rightfully be called shooting by Picard method. A possible application of the new theoretical results is suggested and numerical computer experiments are presented. MATLAB codes are provided.

---

[*] Corresponding Author's Email: faragois@cs.elte.hu

## 1. INTRODUCTION

This chapter considers two-point boundary value problems (TPBVPs) of the form

$$u''(x) = f\big(x, u(x), u'(x)\big), x \in (a, b), \tag{1}$$
$$u(a) = u_a, u(b) = u_b, \tag{2}$$

where $a$, $b$, $u_a$, and $u_b$ are given constants, $u(x)$ is an unknown real-valued function of a real independent variable $x \in [a, b]$, and $f$ is a given function that specifies the differential equation (1) [1-3]. Although there is a lot of literature on TPBVPs, the subject is still an active area of research [4-12]. If the function $f$ is linear with respect to $u$ and $u'$, then the problem is called linear. Otherwise the problem is called nonlinear. In the sequel we use the notation

$$v(x) = u'(x). \tag{3}$$

The following conditions guarantee that the problem (1)-(2) has a solution and that this solution is unique:

$f$ is continuous on the domain $D = \{(x, u, v) \mid x \in [a, b], \ u, v \in \mathbb{R}\}$,

$q = \partial_2 f = \frac{\partial f}{\partial u}$ and $p = \partial_3 f = \frac{\partial f}{\partial v}$ exist and are continuous on $D$,

$q > 0$ on $D$,

$p$ is bounded on $D$.
$$\tag{4}$$

These conditions are sufficient but not necessary. Besides this statement, there are other theorems that guarantee the existence and uniqueness of solution. In the sequel, we assume that the TPBVP (1)-(2) has a unique solution. Often, this solution cannot be defined in a closed form (there is no analytic solution). Therefore, we define a numerical (approximate) solution, which is obtained as the result of some, appropriate for the case, numerical procedure.

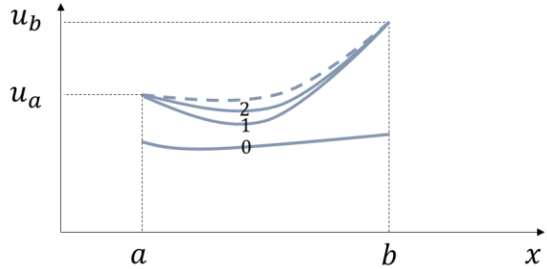

Figure 1. Relaxation method. Three successive relaxation trajectories (functions), marked 0, 1, 2, and the exact solution (dashed line) are shown. The function 0 is arbitrary. The functions 1, 2, … satisfy the two boundary conditions but do not satisfy the differential equation.

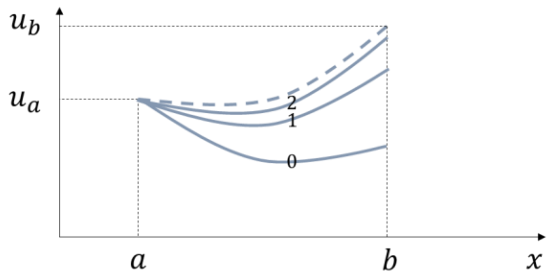

Figure 2. Shooting method. Three successive shooting trajectories (functions), marked 0, 1, 2, and the exact solution (dashed line) are shown. The functions 0, 1, 2, … satisfy the differential equation and the left boundary condition but do not satisfy the right boundary condition.

This work studies numerical methods for solving nonlinear TPBVPs of the form (1)-(2). Two important classes of iterative methods for the solution of such problems are the relaxation methods [13], also known as finite difference methods (FDMs)[1-3,13-17], and the shooting methods [1-3,13,14,18,19]. In the relaxation methods, one starts from an arbitrary function (function 0 in Fig. 1), and each next function (functions 1, 2, … in Fig 1) is brought into a better agreement with the differential equation (1). To achieve this, the differential equation (1) is first discretized using finite differences, and then the resulting nonlinear system, together with the two boundary conditions, is solved by some iterative method for algebraic systems, e.g. Newton method, fixed-point iteration, etc.   Each iterative method generates a particular sequence of functions, called *relaxation trajectories*, that satisfy the two boundary conditions but do not satisfy the differential equation. However, the relaxation trajectories satisfy the differential equation (1) *approximately* and, if the particular iterative method is convergent for the given TPBVP, converge to the exact solution (dashed line in Fig. 1). As shown in section 4, the relaxation trajectories obtained by the FDM can, in principle, be obtained as continuous functions (solid lines in Fig. 1) from the corresponding linearization method. In the shooting methods, the TPBPV (1)-(2) is replaced by an initial value problem (IVP) (Cauchy problem) and a guess is made for the initial condition $v(a)$. The IVP is solved using some numerical technique. The obtained function (function 0 in Fig. 2), called *shooting trajectory*, satisfies the differential equation (1) and the first (left) boundary condition but does not satisfy the second (right) boundary condition. Then, using the end-value of the shooting trajectory and some iterative approach for solving algebraic equations, e.g. Newton method, fixed-point iteration, etc., the initial condition is corrected. The new initial condition is used to obtain a new shooting trajectory (function 1 in Fig. 2), and so on. If the particular iterative method turns out to be convergent for the given TPBVP, then the sequence of shooting trajectories converges to the exact solution (dashed line in Fig. 2).

Other important iterative methods for solving nonlinear TPBVPs are the linearization methods, e.g. the quasilinearization method [20,2,14,17], the Picard linearization method (Picard successive approximations) [21,17], etc. In these methods, one starts from an arbitrary function and uses this function to obtain a linear TPBVP that approximates the nonlinear TPBVP (1)-(2). Then, the linear problem is solved by some technique, and the solution is used to obtain a new linear TPBVP which is a better approximation to (1)-(2). Provided the procedure is convergent, the generated sequence of functions converges to the solution of the nonlinear problem (1)-(2).

In this work we establish some important connections between the linearization methods, the relaxation methods, and the shooting methods. The chapter comprises some of our recent research plus some new unpublished results. The new results are presented in sections 5, 6, and 7.

## 2. ITERATIVE METHODS FOR ALGEBRAIC EQUATIONS

This section considers three important iterative methods for finding roots of nonlinear algebraic equations, namely, the Newton method, the Picard method (fixed-point iteration), and the constant-slope method. The constant-slope method is a particular case of a fixed-point iteration. As shown in the section, all three methods can be viewed as linearization methods because at each iteration step the nonlinear function in the considered nonlinear algebraic equation is replaced by some linear function. The solution to the obtained in this way linear equation is taken as the next root approximation. Consider the algebraic equation

$$F(x) = 0, \tag{5}$$

where $F$ is some given nonlinear function.

Table 2.1. Three different types of linearization (8) for the nonlinear algebraic equation (7).

|  | (i) Newton | (ii) Picard | (iii) Constant-slope |
|---|---|---|---|
| $\rho_k =$ | $\phi'(x_k)$ | 0 | $\phi'(x_0)$ |

Using the nonlinear function $F$, we introduce a new nonlinear function $\phi$:

$$\phi(x) = x - \frac{F(x)}{m}, \tag{6}$$

where $m$ is some fixed number different from zero. Then Eqn. (5) is equivalent to the equation

$$x = \phi(x), \tag{7}$$

i.e. the fixed points of $\phi$ are roots of (5) and vice versa. Let $\tilde{x}$ be a unique root of (5) and let $x_k$ denote a suitable approximation of $\tilde{x}$. The nonlinear function $\phi(x)$ can be approximated by the linear function $\phi(x_k) + \rho_k(x - x_k)$, where $\rho_k$ is some number. We assume that $\rho_k \neq 1$. Thus, instead of equation (7) we write

$$x \approx \phi(x_k) + \rho_k(x - x_k). \tag{8}$$

Clearly, (8) is obtained by linearizing the right-hand side of (7). Depending on the choice of $\rho_k$, we have three different types of linearization: Newton, Picard, or constant-slope (Table 2.1). In the table, $x_0$ denotes some initial guess for $\tilde{x}$. Replacing in (8) the sign $\approx$ by exact equality and, accordingly, $x$ by $x_{k+1}$ gives:

$$x_{k+1} = \phi(x_k) + \rho_k(x_{k+1} - x_k). \tag{9}$$

We can take $x_{k+1}$ as the next approximation to $\tilde{x}$. If iteration (9) turns out to be convergent, then the sequence $x_0, x_1, \ldots$ converges to $\tilde{x}$.

Using (6) in Eqn. (9), we get the usual form of the considered iterative methods:

$$x_{k+1} = x_k - \frac{F(x_k)}{m_k}, \tag{10}$$

where $m_k = m(1 - \rho_k)$. The three different definitions of $\rho_k$ (Table 2.1) yield three different $m_k$-s (Table 2.2) corresponding to the three different iterative methods: Newton, Picard, or constant-slope.

Table 2.2. Three iterative methods (10) for soling the nonlinear equation (5)

|          | (i) Newton | (ii) Picard | (iii) Constant-slope |
|----------|------------|-------------|----------------------|
| $m_k =$  | $F'(x_k)$  | $m$         | $F'(x_0)$            |

As can be seen from (10) and Table 2.2, the constant-slope method is just the Picard method with $m = F'(x_0)$. Both methods are fixed-point iterations and have linear convergence. This means that, for any $k > k_0$ ($k_0$ is a certain integer), $|e_{k+1}| \leq C|e_k|$, where $e_k = x_k - \tilde{x}$ and $0 < C < 1$. The value of $m_k$ for the Newton method is readjusted at each iteration step. It is equal to the slope of the straight line tangent to $F(x)$ at $x = x_k$. The Newton method has quadratic convergence. This means that, if $k > k_0$, $|e_{k+1}| \leq C|e_k|^2$, where $C > 0$. Whether a particular iterative method is convergent or not depends on the nonlinearity $F$ and the initial guess $x_0$ (and the choice of $m$ for the Picard method) [22].

The specific presentation (9) of the Newton, Picard, and constant-slope iterative methods (Table 2.1) is very convenient since it gives us a way to carry out the linearization of the differential equation (1) in an analogous fashion. Note that the right-hand side of (8) for the Newton method is just the first two terms of the Taylor expansion of $\phi(x)$ around $x = x_k$.

## 3. Linearization methods for nonlinear TPBVPs

Let the function $u_k(x), x \in [a, b]$ be some suitable approximation to the solution of the nonlinear problem (1)-(2). Analogously to (8) we can approximate the right-hand side of (1) by a linear function:

$$u''(x) \approx f_k(x) + q_k(x)\big(u(x) - u_k(x)\big) + p_k(x)\big(v(x) - v_k(x)\big), \quad (11)$$

where $v_k(x) = u_k'(x)$, $f_k(x) = f\big(x, u_k(x), v_k(x)\big)$, and $q_k(x), p_k(x)$ are given in Table 3.1. In the table, $u_0(x)$ is some initial guess for the solution of (1)-(2). Adding the boundary conditions (2), replacing in (11) the sign $\approx$ by exact equality and, accordingly, $u(x)$ by its approximation $u_{k+1}(x)$, we get

$$u_{k+1}''(x) = f_k(x) + q_k(x)\big(u_{k+1}(x) - u_k(x)\big) + p_k(x)\big(u_{k+1}'(x) - u_k'(x)\big), \quad (12)$$
$$u_{k+1}(a) = u_a, u_{k+1}(b) = u_b. \quad (13)$$

If the function $u_k(x)$ is given (fixed), then equations (12)-(13) constitute a linear TPBVP for the unknown function $u_{k+1}(x)$.

Table 3.1. Three different types of linearization (11) of the differential equation (1)

|  | (i) Newton (QLM) | (ii) Picard | (iii) Constant-slope |
|---|---|---|---|
| $q_k(x) =$ | $\partial_2 f\big(x, u_k(x), u_k'(x)\big)$ | 0 | $\partial_2 f\big(x, u_0(x), u_0'(x)\big)$ |
| $p_k(x) =$ | $\partial_3 f\big(x, u_k(x), u_k'(x)\big)$ | 0 | $\partial_3 f\big(x, u_0(x), u_0'(x)\big)$ |

Thus, if we start from some initial guess $u_0(x)$, we can solve successively (12)-(13) for $k = 0,1, \ldots$ and obtain $u_1(x), u_2(x), \ldots$ . If the sequence of solutions $u_1(x), u_2(x), \ldots$ is convergent, it converges to the solution of the nonlinear problem (1)-(2). The method (12)-(13) for the Newton case (column (i) of Table 3.1) is the well-known quasilinearization method (QLM). Originally, it was introduced by Bellman and Kalaba [20] as a generalization of the Newton-Raphson method. It can be viewed as Newton method on operator level [2]. We note that, under conditions (4) and provided that $u_k(x)$ is continuously differentiable, the linear problem (12)-(13) arising from quasilinearization has a unique solution. The Picard linearization method and the constant-slope linearization method are also given by equations (12)-(13) but with $q_k(x)$, $p_k(x)$ given in columns (ii) and (iii)  of Table 3.1.

To solve the linear problem (12)-(13) we can apply different numerical techniques. In the next section we prove that applying the finite difference method leads to the usual iteration formula of the respective relaxation method. Then we show that the shooting-projection iteration, with (12)-(13) as a projection transformation, leads to the respective shooting method. Thus, the linearization methods (12)-(13) can be used as a basis to derive both the relaxation and the shooting methods. A possible application of these theoretical results is suggested in the final sections.

## 4. DERIVING THE RELAXATION METHODS

Let us divide the interval $[a, b]$ by $N$ equally separated mesh-points:

$$x_i = a + (i - 1)h, h = \frac{b - a}{N - 1}, i = 1,2, \ldots, N, \tag{14}$$

where $N > 1$. The points $x_i$ (14), where $x_1 = a$ and $x_N = b$, define a uniform mesh on the interval $[a, b]$. We are going to prove that using finite differences to discretize (12)-(13) on the mesh (14) yields the usual iteration formula of the respective relaxation method. Let $u_i^k$ and $u_i^{k+1}$ denote, respectively, approximations of the values of $u_k(x)$ and $u_{k+1}(x)$ at $x = x_i$ (14).

Table 4.1. Three different relaxation methods (25) for solving the nonlinear TPBVP (1)-(2)

|  | (i) Newton | (ii) Picard | (iii) Constant-slope |
|---|---|---|---|
| $q_i^k =$ | $\partial_2 f\big(x_i, u_i^k, \mathcal{D}u_i^k\big)$ | 0 | $\partial_2 f\big(x_i, u_i^0, \mathcal{D}u_i^0\big)$ |
| $p_i^k =$ | $\partial_3 f\big(x_i, u_i^k, \mathcal{D}u_i^k\big)$ | 0 | $\partial_3 f\big(x_i, u_i^0, \mathcal{D}u_i^0\big)$ |

To approximate the first derivatives of $u_k(x)$ and $u_{k+1}(x)$ at the inner mesh-points, we could use forward (a), backward (b), or central (c) difference approximation:

$$\text{(a)} \quad \mathcal{D}_+ u_i^k = \frac{u_{i+1}^k - u_i^k}{h}, \mathcal{D}_+ u_i^{k+1} = \frac{u_{i+1}^{k+1} - u_i^{k+1}}{h},$$

$$\text{(b)} \quad \mathcal{D}_- u_i^k = \frac{u_i^k - u_{i-1}^k}{h}, \mathcal{D}_- u_i^{k+1} = \frac{u_i^{k+1} - u_{i-1}^{k+1}}{h}, \quad (15)$$

$$\text{(c)} \quad \mathcal{D}_0 u_i^k = \frac{u_{i+1}^k - u_{i-1}^k}{2h}, \mathcal{D}_0 u_i^{k+1} = \frac{u_{i+1}^{k+1} - u_{i-1}^{k+1}}{2h}.$$

Discretizing (12) using the central difference approximation to approximate the second derivative of $u_{k+1}(x)$, we get

$$\frac{u_{i-1}^{k+1} - 2u_i^{k+1} + u_{i+1}^{k+1}}{h^2} = f_i^k + q_i^k\big(u_i^{k+1} - u_i^k\big) + p_i^k\big(\mathcal{D}u_i^{k+1} - \mathcal{D}u_i^k\big),$$

$$i = 2, 3, \ldots, N-1, \quad (16)$$

where $\mathcal{D} \in \{\mathcal{D}_+, \mathcal{D}_-, \mathcal{D}_0\}$, $f_i^k = f\big(x_i, u_i^k, \mathcal{D}u_i^k\big)$, and $q_i^k, p_i^k$ are given in Table 4.1. Note that, since $h$ is a constant, $\mathcal{D}u_i^k$ and $\mathcal{D}u_i^{k+1}$ are linear homogenous functions of their arguments. Therefore, according to Euler's theorem on homogenous functions,

$$\mathcal{D}u_i^k = \frac{\partial \mathcal{D}u_i^k}{\partial u_{i-1}^k} u_{i-1}^k + \frac{\partial \mathcal{D}u_i^k}{\partial u_i^k} u_i^k + \frac{\partial \mathcal{D}u_i^k}{\partial u_{i+1}^k} u_{i+1}^k, \quad (17)$$

$$\mathcal{D}u_i^{k+1} = \frac{\partial \mathcal{D}u_i^{k+1}}{\partial u_{i-1}^{k+1}} u_{i-1}^{k+1} + \frac{\partial \mathcal{D}u_i^{k+1}}{\partial u_i^{k+1}} u_i^{k+1} + \frac{\partial \mathcal{D}u_i^{k+1}}{\partial u_{i+1}^{k+1}} u_{i+1}^{k+1}. \quad (18)$$

Substituting (17) and (18) into (16) and using the obvious relations

$$\frac{\partial \mathcal{D}u_i^{k+1}}{\partial u_{i-1}^{k+1}} = \frac{\partial \mathcal{D}u_i^k}{\partial u_{i-1}^k}, \frac{\partial \mathcal{D}u_i^{k+1}}{\partial u_i^{k+1}} = \frac{\partial \mathcal{D}u_i^k}{\partial u_i^k}, \frac{\partial \mathcal{D}u_i^{k+1}}{\partial u_{i+1}^{k+1}} = \frac{\partial \mathcal{D}u_i^k}{\partial u_{i+1}^k}, \quad (19)$$

equation (16) is transformed into the form

$$L_{i,i-1}^k u_{i-1}^{k+1} + L_{i,i}^k u_i^{k+1} + L_{i,i+1}^k u_{i+1}^{k+1} = L_{i,i-1}^k u_{i-1}^k + L_{i,i}^k u_i^k + L_{i,i+1}^k u_{i+1}^k - G_i^k, \quad (20)$$

where

$$\left. \begin{aligned} L_{i,i-1}^k &= \left(1 - h^2 p_i^k \frac{\partial \mathcal{D}u_i^k}{\partial u_{i-1}^k}\right), \\ L_{i,i}^k &= \left(-2 - h^2 q_i^k - h^2 p_i^k \frac{\partial \mathcal{D}u_i^k}{\partial u_i^k}\right), \\ L_{i,i+1}^k &= \left(1 - h^2 p_i^k \frac{\partial \mathcal{D}u_i^k}{\partial u_{i+1}^k}\right), \end{aligned} \right\} \quad (21)$$

and

$$G_i^k = u_{i-1}^k - 2u_i^k + u_{i+1}^k - h^2 f_i^k. \quad (22)$$

Equation (20), together with the boundary conditions

$$u_1^{k+1} = u_a, u_N^{k+1} = u_b, \quad (23)$$

can be written in a matrix form

$$\mathbf{L}^k \mathbf{u}^{k+1} = \mathbf{L}^k \mathbf{u}^k - \mathbf{G}^k. \quad (24)$$

The components of the column-vector $\mathbf{G}^k$ are $G_1^k = u_1^k - u_a$, $G_N^k = u_N^k - u_b$, and (22) for $i = 2,3,\ldots,N-1$. The nonzero components of the $N \times N$ matrix $\mathbf{L}^k$ are $L_{1,1}^k = 1$, $L_{N,N}^k = 1$, and (21) for $i = 2,3,\ldots,N-1$. In (24) $\mathbf{u}^k = \left[u_1^k, u_2^k, \ldots, u_N^k\right]^T$ and $\mathbf{u}^{k+1} = \left[u_1^{k+1}, u_2^{k+1}, \ldots, u_N^{k+1}\right]^T$. Under conditions (4) and providing $h < 2/P$ when $\mathcal{D} = \mathcal{D}_0$ or $h < 1/P$ when $\mathcal{D} \in \{\mathcal{D}_+, \mathcal{D}_-\}$, where $P$ is the upper bound of $|p|$, the matrix $\mathbf{L}^k$ is strictly diagonally dominant, hence nonsingular. Multiplying (24) from the left by the inverse of $\mathbf{L}^k$, we finally get

$$\mathbf{u}^{k+1} = \mathbf{u}^k - (\mathbf{L}^k)^{-1}\mathbf{G}^k. \quad (25)$$

Formula (25) is nothing else but the iteration formula of the respective relaxation method. When (12)-(13) is the Newton linearization method (quasilinearization), then (25) is the Newton relaxation method (see columns (i) in Tables 3.1 and 4.1). When (12)-(13) is the Picard linearization method, then (25) is the Picard relaxation method (see columns (ii) in Tables 3.1 and 4.1). When (12)-(13) is the constant-slope linearization method, then (25) is the constant-slope relaxation method (see columns (iii) in Tables 3.1 and 4.1). The usual way to derive formula (25) is first to discretize the nonlinear problem (1)-(2) using the FDM and then to apply one of three considered iterative methods to the obtained nonlinear algebraic system. As we have just proven, the relaxation methods (25) can be derived directly from the linearization methods (12)-(13) by simply discretizing equations (12)-(13) with the FDM. Therefore, given a function $u_k$, we can obtain $u_{k+1}$ by either applying (25) or solving (12)-(13) by some technique other than the FDM. No matter which way we use, the resulting function $u_{k+1}$, called *relaxation trajectory*, is, within numerical precision, the same. This equivalence was first proven by the authors in [17]. The result was used to construct a way of replacing the relaxation methods (25) by successive application of the linear shooting method.

## 5. DERIVING THE SHOOTING METHODS

This and the following sections of the chapter present some unpublished new results of our recent research work. In this section we show that the linearization methods (12)-(13) can be used a basis to derive the shooting methods. First, we demonstrate that if the function $u_k$ is a shooting trajectory (IVP solution), then equations (12)-(13) define a projection transformation of $u_k$ into $u_{k+1}$. In this case the relaxation trajectory $u_{k+1}$ is called a projection trajectory. Then, we introduce the shooting-projection iteration. In essence, the shooting-projection iteration is the following: the shooting trajectory is transformed into a projection trajectory, then the first derivative of the projection trajectory at the left boundary is used as a new initial condition and a new shooting trajectory is found, and so on. We show that applying the shooting-projection iteration with the quasi-linearization equation as a projection transformation results in the usual shooting by Newton method. When the Picard linearization is used as a projection, then we get the recently proposed by the authors shooting-projection method [19], while the constant-slope linearization leads to the shooting by constant-slope method.

### 5.1. Relaxing the shooting trajectory – projection transformation

Let $u(x; v_a^k), x \in [a, b], v_a^k \in \mathbb{R}$ be a solution to the following initial value problem:

$$u''(x; v_a^k) = f\left(x, u(x; v_a^k), u'(x; v_a^k)\right), \tag{26}$$
$$u(a; v_a^k) = u_a, u'(a; v_a^k) = v_a^k. \tag{27}$$

The function $u(x; v_a^k)$ is called a shooting trajectory. It satisfies the differential equation (1) and the left boundary condition in (2) but, typically, does not satisfy the right boundary condition in (2). The number $v_a^k$ is some approximation for the first derivative at $x = a$ of the solution of the nonlinear problem (1)-(2). Let us relax the shooting trajectory into $u_{k+1}(x)$ using the linearization equations (12)-(13). In other words, we replace the function $u_k(x)$ in (12)-(13) by the shooting-trajectory $u(x; v_a^k)$. Taking into account (26), equations (12)-(13) become

$$u''_{k+1}(x) = u''_k(x) + q_k(x)\big(u_{k+1}(x) - u_k(x)\big) + p_k(x)\big(u'_{k+1}(x) - u'_k(x)\big), \text{ (28)}$$
$$u_{k+1}(a) = u_a, u_{k+1}(b) = u_b, \tag{29}$$

where $u_k(x) = u(x; v_a^k)$. Let $\hat{P}$ be an operator defined by the equation

$$\hat{P}u_k = u_{k+1}, \tag{30}$$

where $u_k$ is *any* fixed function (not necessarily a shooting trajectory) for which (28)-(29) has a unique solution and $u_{k+1}$ is this unique solution. Obviously, if $u_k$ satisfies the boundary conditions, then the unique solution to (28)-(29) is just $u_{k+1} = u_k$, i.e. $\hat{P}u_k = u_k$. Since $u_{k+1}$ satisfies the boundary conditions, it follows that

$$\hat{P}u_{k+1} = u_{k+1}. \tag{31}$$

Hence, from (30) and (31) we get

$$\hat{P}^2 u_k = \hat{P}u_k. \tag{32}$$

Since this is true for any $u_k$, it follows that $\hat{P}^2 = \hat{P}$, i.e. $\hat{P}$ is idempotent. In general, an idempotent mapping of a set (or other mathematical structure) into a subset (or a sub-structure) is a projection [23]. Therefore, the transformation of $u_k$ into $u_{k+1}$ defined by (28)-(29) can be considered a projection transformation. [1]

## 5.2. Shooting-projection iterative procedure

As shown in the previous section, the linear problem (28)-(29) defines a projection transformation of $u_k$ into $u_{k+1}$. Depending on what values of $q_k(x)$, $p_k(x)$ we choose (Table 3.1), we have three types of projection: Newton, Picard, or constant-slope. The function $u_{k+1}$, i.e. the solution to (28)-(29), is called a projection trajectory. The projection trajectory satisfies the boundary conditions (2) exactly and the differential equation (1) approximately. Therefore, it is an approximate solution to the nonlinear TPBVP (1)-(2). Hence, we could use the first derivative of the projection trajectory at the first (left) boundary as a new initial condition:

$$v_a^{k+1} = u'_{k+1}(a). \tag{33}$$

---

[1] Strictly speaking, in linear algebra and functional analysis, an operator must be linear and idempotent to be a projection. The operator $\hat{P}$ is linear if the boundary conditions are homogenous. However, all the results and conclusions reached in this chapter hold for the non-homogenous BCs (2) as well.

Then, we could use the new initial condition $v_a^{k+1}$ to find a new shooting trajectory, and so on. This is the shooting-projection iterative procedure or *shooting-projection iteration* (SPI). In summary, the algorithm of the SPI is:

1) Use $v_a^k$ and find the shooting trajectory $u(x; v_a^k)$ that satisfies (26)-(27).
2) Solve (28)-(29) with $u_k(x) = u(x; v_a^k)$ to find the projection trajectory $u_{k+1}$.
3) Use (33) to find the new initial condition $v_a^{k+1}$ and repeat 1)-3) with $v_a^{k+1}$.

If the iteration is convergent, then $u(x; \tilde{v}_a)$, where $\tilde{v}_a = \lim v_a^k$ as $k \to \infty$, is the sought solution to the nonlinear TPBVP (1)-(2).

According to the results of section 4, to find the projection trajectory, i.e. to solve (28)-(29), we can apply (25) in which $\mathbf{u}^k$ represents the shooting trajectory $u(x; v_a^k)$. We note that, since the shooting trajectory satisfies the first boundary condition, i.e. $u_1^k = u_a$, and the nonlinear differential equation (1), i.e. in a discretized form $u_{i-1}^k - 2u_i^k + u_{i+1}^k = h^2 f_i^k$, as $\mathbf{G}^k$ in (25) we can take

$$\mathbf{G}^k = [0, 0, \ldots, 0, u(b; v_a^k) - u_b]^T. \tag{34}$$

The option of calculating the projection trajectory via (25) will be used later on in the chapter when we consider possible applications of the theoretical results.

It turns out that, for all three considered SPI methods, namely Newton, Picard, and constant-slope, it is possible to find a relation between the $k$-th and the improved $(k+1)$-th guess for the initial derivative values, namely $v_a^{k+1}$ as a function of $v_a^k$. These relations are derived in the next three sections.

## 5.3.  Deriving the shooting by Newton method

Let

$$y(x) = u_{k+1}(x) - u_k(x), \tag{35}$$

where $u_k(x) = u(x; v_a^k)$ is a shooting trajectory and $u_{k+1}(x)$ is the corresponding Newton projection trajectory satisfying (28)-(29). Using (35), equations (28)-(29) can be written as

$$y''(x) = q_k(x)y(x) + p_k(x)y'(x), \tag{36}$$
$$y(a) = 0, y(b) = u_b - u(b; v_a^k). \tag{37}$$

where $q_k(x), p_k(x)$ are given in Table 3.1 column (i). Taking the derivative of (35) and using (33) and the second initial condition in (27), we replace the boundary conditions (37) with the initial conditions

$$y(a) = 0, y'(a) = v_a^{k+1} - v_a^k. \tag{38}$$

Obviously, the initial value problem (36), (38) is equivalent to the two-point boundary value problem (36), (37). Now, introducing the function $z(x)$ such that

$$y(x) = (v_a^{k+1} - v_a^k)z(x), \tag{39}$$

the initial value problem (36), (38) becomes

$$z''(x) = q_k(x)z(x) + p_k(x)z'(x), \tag{40}$$
$$z(a) = 0, z'(a) = 1. \tag{41}$$

Differentiating (26)-(27) with respect to $v_a^k$ (as in [18]), comparing the result to (40)-(41), and taking into account that (40)-(41) has a unique solution, we get

$$z(x) = \frac{\partial u(x; v_a^k)}{\partial v_a^k}. \tag{42}$$

At $x = b$, using (39) and the second boundary condition in (37), equation (42) yields

$$v_a^{k+1} = v_a^k - \frac{u(b; v_a^k) - u_b}{\frac{\partial u(b; v_a^k)}{\partial v_a^k}}. \tag{43}$$

Equation (43) is the well-known iteration formula of the shooting by Newton method. Therefore, the shooting-projection iteration with projection transformation (28)-(29), where $q_k(x), p_k(x)$ are given in Table 3.1 column (i),  is equivalent to the shooting by Newton method.

### 5.4.   Deriving the shooting by Picard method

The Picard projection transformation (28)-(29) with $q_k(x)$, $p_k(x)$ given in Table 3.1 column (ii) is

$$u''_{k+1}(x) = u''_k(x), \tag{44}$$

together with the boundary conditions (29), where $u_k(x) = u(x; v^k_a)$ is a shooting trajectory and $u_{k+1}(x)$ is the corresponding projection trajectory. Integrating (44) on $[a, x]$, and then integrating the result on $[a, b]$, we get

$$u_{k+1}(b) - u_{k+1}(a) - u'_{k+1}(a)(b-a) = u_k(b) - u_k(a) - u'_k(a)(b-a). \tag{45}$$

Taking into account that the shooting trajectory $u_k(x) = u(x; v^k_a)$ satisfies the initial conditions (27) and the projection trajectory $u_{k+1}(x)$ satisfies the boundary conditions (29), we get

$$v^{k+1}_a = v^k_a - \frac{u(b; v^k_a) - u_b}{b-a}. \tag{46}$$

Formula (46) is the iteration formula of the recently proposed shooting-projection method [19]. Since the projection transformation (44), (29) stems from the Picard linearization method (12)-(13) Table 3.1 column (ii), the shooting-projection method (46) can rightfully be called *shooting by Picard method*.

### 5.5.   Deriving the shooting by constant-slope method

Using again the functions $y(x)$ (35) and $z(x)$ defined by (39), the initial value problem (40)-(41) for the constant-slope case is derived:

$$z''(x) = q_0(x)z(x) + p_0(x)z'(x), \tag{47}$$
$$z(a) = 0, z'(a) = 1, \tag{48}$$

where $q_0(x)$, $p_0(x)$ are given in Table 3.1 column (iii). Differentiating (26)-(27) for $k = 0$ with respect to $v^0_a$, comparing the result to (47)-(48), and taking into account that the solution to (47)-(48) is unique, we get

$$z(x) = \frac{\partial u(x; v^0_a)}{\partial v^0_a}. \tag{49}$$

Table 5.2.1. Three different iteration formulas (51) corresponding to the three different projections (28)-(29)

|         | (i) Newton | (ii) Picard | (iii) Constant-slope |
|---------|------------|-------------|----------------------|
| $l_k =$ | $\dfrac{\partial u(b; v_a^k)}{\partial v_a^k}$ | $b - a$ | $\dfrac{\partial u(b; v_a^0)}{\partial v_a^0}$ |

Finally, using (39) and the second boundary condition in (37), equation (49), at $x = b$, yields

$$v_a^{k+1} = v_a^k - \frac{u(b; v_a^k) - u_b}{\dfrac{\partial u(b; v_a^0)}{\partial v_a^0}}. \tag{50}$$

Equation (50) is the iteration formula of the shooting by constant-slope method. Therefore, the shooting-projection iteration with projection transformation (28)-(29), where $q_k(x)$, $p_k(x)$ are given in Table 3.1 column (iii),  is equivalent to the shooting by constant-slope method.

## 5.6.   General iteration formula of the SPI

Formulas (43), (46), and (50) can be unified in one single formula:

$$v_a^{k+1} = v_a^k - \frac{u(b; v_a^k) - u_b}{l_k}, \tag{51}$$

where $l_k$ depends on the choice of projection: Newton, Picard, or constant-slope (Table 5.2.1). This is the shooting-projection iteration formula for correcting the initial condition. It allows us to avoid (skip) the explicit calculation of the projection trajectory. As discussed, iteration (51), with $l_k$ given in Table 5.2.1, is the same as the iteration formulas of known shooting methods, namely: (i) shooting by Newton method, (ii) shooting-projection method [19], and (iii) shooting by constant-slope method. The theoretical importance of the result is that it shows that the SPI with the usual linearization methods as projections leads to the respective shooting methods. The practical importance of the result is that it gives us an alternative of calculating $v_a^{k+1}$ not through formula (51) but via the projection trajectory. This is discussed in detail in the next sections.

## 6. APPLICATION OF THE THEORETICAL RESULTS

Let $u(x; v_a^k)$ be the shooting trajectory satisfying the initial value problem (26)-(27). Using the iteration formula of the shooting by Newton method (43) we can obtain a new initial condition $v_a^{k+1}$. According to the result in section 5.3, the value of $v_a^{k+1}$ is exactly equal to the first derivative at $x = a$ of the corresponding projection trajectory $u_{k+1}(x)$ obtained from the projection transformation (28)-(29) with $q_k(x)$, $p_k(x)$ given in Table 3.1 column (i). As discussed in section 5.1, this projection transformation is just the quasilineartization (12)-(13) Table 3.1 column (i) applied on the shooting trajectory. This result gives us an alternative to finding the new initial condition $v_a^{k+1}$, namely, instead of using (43), we could use (33). In order to use (33) we need to find the projection trajectory $u_{k+1}(x)$ that satisfies (28)-(29) . According to the result in section 4, the projection trajectory can be obtained in a discretized form, i.e. $\mathbf{u}^{k+1}$, by applying the Newton relaxation equation (25), where $q_i^k$, $p_i^k$ are given in Table 4.1 column (i) and $\mathbf{G}^k$ is given by (34). Then, we can use

$$\frac{u_2^{k+1} - u_1^{k+1}}{h} \tag{52}$$

as an approximation to $v_a^{k+1}$.

This section discusses a possible situation in which the option of calculating $v_a^{k+1}$ via (52) could be very useful. Suppose we need to solve a TPBVP of the form

$$u''(x) = f\big(u(x)\big), \tag{53}$$

$$u(a) = u_a, u(b) = u_b, \tag{54}$$

but for some reason the derivative of $f$ with respect to $u$ is not available. For example, if $u$ is a position and $x$ is time, then (53) is the equation of motion for a particle (small object) of unit mass travelling under the influence of the force $f$. Launching the particle with initial velocity $v_a^k$ and tracing its motion is equivalent to finding the shooting trajectory $u(x; v_a^k)$ , i.e. to solving the equation of motion (53) with initial conditions (27). If the particle does not arrive at the prescribed position $u_b$ at time $b$, which is typically the case, we are faced with the need to correct (change) the initial condition $v_a^k$ to a new one $v_a^{k+1}$. However, if we want to use the shooting by Newton method (43) we need to solve the IVP (40)-(41) and obtain $\partial u(b; v_a^k)/\partial v_a^k$ needed for formula (43) from $z(x)$ at $x = b$ (see Eqn. (42)). Since $q = \partial f/\partial u$ is not available, this is impossible (at least not in a direct way). Another (approximate) method would be to launch a second trial with initial velocity $v_a^k + \Delta v_a^k$, where $\Delta v_a^k$ is a small change in the initial velocity,  and then calculate $\Delta u(b; v_a^k) / \Delta v_a^k$. However, if launching a second trial is impossible or undesirable, then the described in this paragraph new approach may come in very handy. Taking the derivative of (53) we get

$$u'''(x) = \partial_u f\big(u(x)\big) u'(x), \tag{55}$$

i.e., using $v$ (3) and $q$,

$$v''(x) = q(x) v(x). \tag{56}$$

For any trajectory $u(x)$ that satisfies the differential equation (53), equation (56) must hold. Discretizing equation (56) on the mesh (14) using the central difference approximation for $v''(x)$, we get

$$\frac{v_{i-1} - 2v_i + v_{i+1}}{h^2} = q_i v_i, i = 2,3, \dots, N-1, \tag{57}$$

where $q_i = \partial_u f(u_i)$. Let $\mathbf{u}^k = \big[u_1^k, u_2^k, \dots, u_N^k\big]^T$ represent the shooting trajectory $u(x; v_a^k)$ in a discretized form, and let $\mathbf{v}^k = \big[v_1^k, v_2^k, \dots, v_N^k\big]^T$ represent its first derivative $v(x; v_a^k) = u'(x; v_a^k)$. Since the shooting trajectory satisfies the differential equation (53), it follows that $v_i^k$ satisfy (57). Hence, after substituting $v_i^k$ into (57) and rearranging, we get

$$-2 - h^2 q_i^k = -\frac{v_{i-1}^k + v_{i+1}^k}{v_i^k}, i = 2,3, \dots, N-1, \tag{58}$$

where $q_i^k = \partial_u f(u_i^k)$. For the TPBVP (53)-(54), since $p = \partial_v f = 0$, the nonzero elements of the Jacobian matrix $\mathbf{L}^k$ for the Newton relaxation method are

$$L_{1,1}^k = 1, L_{N,N}^k = 1,$$
$$L_{i,i-1}^k = 1, L_{i,i}^k = -2 - h^2 q_i^k, L_{i,i+1}^k = 1, i = 2,3, \dots, N-1. \tag{59}$$

Hence, using (58), the only unknown elements of the matrix $\mathbf{L}^k$ can be found from the shooting trajectory itself:

$$L_{i,i}^k = -\frac{v_{i-1}^k + v_{i+1}^k}{v_i^k}, i = 2,3, \dots, N-1. \tag{60}$$

Table 7.1. Solving the TPBVP (61)-(62) by the traditional shooting by Newton method (43)

| $k$ | $v_a^k$ – initial condition | $E_k$ – deviation from 2$^{nd}$ BC |
|---|---|---|
| 0 | 0 | -0.433349035739307 |
| 1 | 0.379948530223661 | 0.026009489270876 |
| 2 | 0.359783026933729 | 0.000100006717963 |

Table 7.2. Solving the TPBVP (61)-(62) by the proposed new approach (52)

| $k$ | $v_a^k$ – initial condition | $E_k$ – deviation from 2$^{nd}$ BC |
|---|---|---|
| 0 | 0 | -0.433349035739307 |
| 1 | 0.379942276669709 | 0.026001423439514 |
| 2 | 0.359786457564626 | 0.000104397665347 |

Of course, equation (60) can be used when $v_i^k \neq 0$. If $v_i^k = 0$ for some particular $i$, we could simply extrapolate the value of $L_{i,i}^k$ from the values at the neighboring mesh-points. In the example provided in the next section, there are no such complications. Having calculated the matrix elements (60), we could use the Newton relaxation (25) with (34) for $\mathbf{G}^k$ to transform the shooting trajectory $\mathbf{u}^k$ into the projection trajectory $\mathbf{u}^{k+1}$. Then, (52) gives the next initial condition $v_a^{k+1}$ for the shooting by Newton method.

## 7. Computer experiments and Matlab codes

Consider the TPBVP

$$u''(x) = \left(u(x)\right)^3, \tag{61}$$

$$u(0) = 1/2, u(1) = 1. \tag{62}$$

The problem is solved by the shooting by Newton method, first by using the traditional approach to correct the initial condition, i.e. Eqn. (43), and then by the new approach described in the previous section, namely (52) to approximate $v_a^{k+1}$. For starting initial condition we choose $v_a^0 = 0$. To find the shooting trajectory, at each iteration step, the explicit Euler method EE_ [17] is used. The discretization step is $h = 0.001$. The deviation from the second boundary condition is denoted by

$$E_k = u(b; v_a^k) - u_b. \tag{63}$$

The results are shown in Table 7.1 for the traditional shooting by Newton method and in Table 7.2 for the new approach.

Table 7.3. MATLAB codes for the solution of the TPBVP (61)-(62)

| Traditional shooting by Newton (43) | The proposed new approach (52) |
|---|---|

```
function main
N=1001;
a=0;  b=1;
ua=0.5; ub=1;
va=0;
h=(b-a)/(N-1);
x=zeros(N,1); u=zeros(N,1);
v=zeros(N,1);
z=zeros(N,1); w=zeros(N,1);
for i=1:N
 x(i)=a+h*(i-1);
end
u(1)=ua; v(1)=va;
z(1)=0; w(1)=1;
hold on; format long;
E=1;
while(abs(E)>0.001)
 for i=2:N
  u(i)=u(i-1)+h*v(i-1);
  v(i)=v(i-1)+h*f(x(i),u(i),v(i-1));
  z(i)=z(i-1)+h*w(i-1);
  qi=q(x(i),u(i),v(i-1));
  pi=p(x(i),u(i),v(i-1));
  w(i)=w(i-1)+h*(qi*z(i)+pi*w(i-1));
 end
 plot(x,u);
 E=u(N)-ub
 v(1)=v(1)-E/z(N);
 end
end

function f_=f(x,u,v)
 f_=u*u*u;
end

function p_=p(x,u,v)
 p_=0;
end

function q_=q(x,u,v)
 q_=3*u*u;
end
```

```
function main
N=1001;
a=0;  b=1;
ua=0.5; ub=1;
va=0;
h=(b-a)/(N-1);
x=zeros(N,1); u=zeros(N,1);
v=zeros(N,1);
L=zeros(N,N); G=zeros(N,1);
for i=1:N
 x(i)=a+h*(i-1);
end
L(1,1)=1; L(N,N)=1;
for i=2:N-1
 L(i,i-1)=1; L(i,i+1)=1;
end
u(1)=ua; v(1)=va;
hold on; format long;
E=1;
while(abs(E)>0.001)
 for i=2:N
  u(i)=u(i-1)+h*v(i-1);
  v(i)=v(i-1)+h*f(x(i),u(i),v(i-1));
 end
 for i=2:N-1
  L(i,i)=-(v(i-1)+v(i+1))/v(i);
 end
 plot(x,u,'-g');
 E=u(N)-ub
 G(N)=E;
 uProj=u-L\G;
 v(1)=(uProj(2)-uProj(1))/h;
 end
end

function f_=f(x,u,v)
 f_=u*u*u;
end
```

MATLAB codes for the solution of (61)-(62) by the traditional shooting by Newton method and by the proposed new approach are presented in Table 7.3. In the first code we have used the notation $w = z'$. Note that, for the proposed new approach, the derivatives $q = \partial_u f$ and $p = \partial_v f$ are not defined at all.

## 8. Conclusion

This chapter revealed important connections between the linearization, the relaxation, and the shooting methods for nonlinear TPBVPs. It was demonstrated that discretizing the linearization equations by using the FDM yields the respective relaxation method. It was also shown that the shooting-projection iteration with the linearization equations as projections leads to the usual shooting methods. An application of the theoretical results was proposed whereby the shooting by Newton method is carried out not through the traditional iteration formula but via explicit calculation of the projection trajectory. The approach is useful because for certain TPBVPs it requires knowledge only of the nonlinear function $f$ but not of its partial derivatives.

### Acknowledgement

This research was partially carried out in the ELTE Institutional Excellence Program (1783-3/2018/FEKUTSRAT) supported by the Hungarian Ministry of Human Capacities, and it was supported by the Hungarian Scientific Research Fund SNN125119. The research was also supported by the program for young scientists and postdocs in Bulgaria.